\documentclass[a4paper,10pt]{article}

\usepackage{amssymb, amsmath, amsthm, latexsym, verbatim, enumerate,array}

\renewcommand \a{\alpha}

\newcommand \K{\delta}

\newcommand \la{\lambda}
\newcommand \id{\mathrm{id}}
\newcommand \Rn{\mathbb R^n}

\newcommand \Span{\mathrm{Span}}
\newcommand \Tr{\mathrm{Tr} \, }

\newcommand \<{\langle}
\renewcommand \>{\rangle}

\newcommand \g{\mathfrak{g}}
\newcommand \p{\mathfrak{p}}
\newcommand \f{\mathfrak{f}}
\newcommand \n{\mathfrak{n}}

\newcommand \ag{\mathfrak{a}}
\newcommand \ad{\mathrm{ad}}
\newcommand \ric{\mathrm{ric}}
\newcommand \Ric{\mathrm{Ric}}
\newcommand \Der{\mathrm{Der}}
\newcommand \End{\mathrm{End}}

\newcommand \GL{\mathrm{GL}}

\theoremstyle{plane}

\newtheorem*{theorem*}{Theorem}
\newtheorem*{corollary*}{Corollary}
\newtheorem{lemma}{Lemma}

\newtheorem*{proposition*}{Proposition}

\newtheorem*{namedtheorem}{\theoremname}
\newcommand{\theoremname}{te}

\theoremstyle{definition}

\newtheorem*{definition*}{Definition}

\theoremstyle{remark}
\newtheorem{remark}{Remark}
\newtheorem*{example}{Example}

\newcommand\cn{\centering}

\begin{document}

\title{Einstein solvmanifolds with free nilradical}

\author{Y.Nikolayevsky}

\date{}
%\keywords{}
%subjclass[2000]Primary: 53C25
\maketitle

\begin{abstract}
We classify solvable Lie groups with free nilradical admitting an Einstein left-invariant metric. Any such group is
essentially determined by the nilradical of its Lie algebra, which is then called an Einstein nilradical. We show that
among the free Lie algebras, there are very few Einstein nilradicals. Except for the one-step (abelian) and the
two-step ones, there are only six others: $\f(2,3), \; \f(2,4), \; \f(2,5), \; \f(3,3)$, $\f(4,3), \; \f(5,3)$ (here
$\f(m,p)$ is a free $p$-step Lie algebra on $m$ generators). The reason for that is the inequality-type restrictions on
the eigenvalue type of an Einstein nilradical obtained in the paper.
\end{abstract}

\section{Introduction}
\label{s:intro}

The theory of Riemannian homogeneous spaces with an Einstein metric splits into three very different cases
depending on the sign of the Einstein constant, the scalar curvature. Among them, only in the Ricci-flat case
the picture is complete: by the result of \cite{AK}, every Ricci-flat homogeneous space is flat.
%(add more refs on recognition of such a space? BB? also ref on viberg/alekseevski/soldovnikov from besse).
Despite many powerful existence, nonexistence and classification results (see \cite{WZ, BWZ}), the positive scalar
curvature case is far from being understood. In the case of negative scalar curvature (which we are dealing with in
this paper), one has two major open questions. The first one is the \emph{Alekseevski Conjecture} \cite{A1} asserting
that an Einstein homogeneous Riemannian space with negative scalar curvature admits a simply transitive solvable
isometry group. This is equivalent to saying that any such space is a \emph{solvmanifold}, a solvable Lie group with a
left-invariant Riemannian metric satisfying the Einstein condition.

The second question arises from the fact that all the known examples of Einstein solvmanifolds are \emph{standard}.
This means that the metric solvable Lie algebra $\g$ of such a solvmanifold has the following property: the orthogonal
complement $\ag$ to the derived algebra of $\g$ is abelian. \emph{Are there any nonstandard Einstein solvmanifolds?}
In the paper \cite{H}, together with a deep and detailed analysis of Einstein solvmanifolds, it was
shown that for several classes of solvmanifolds, the existence of an Einstein left invariant metric implies standardness.
What is more, nonstandard solvmanifolds must be algebraically very different from the standard ones (they can not even
share the same nilradical \cite{Ni}). Further progress was made in \cite{S}, for solvmanifolds with
$\dim[\ag, \ag] = 1$. Einstein solvmanifolds of dimension up to $6$ are classified in \cite{N1, NN} (all of them are
standard). In \cite{Ni} it is proved that there are no nonstandard solvmanifolds of dimension less than ten.

The problem of studying and classifying standard Einstein solvmanifolds (standard metric solvable Einstein Lie algebras)
was given a great deal of attention in the last decade \cite{F, GK, L1, L2, L3, LW, N1, P, T}.

Any such solvmanifold admits a rank-one reduction \cite{H}. On the Lie algebra
level, this means that if $\g$ is a standard Einstein metric solvable Lie algebra and $\g=\ag \oplus \n$
(orthogonal sum), with $\n$ the nilradical of $\g$, then there exists a one-dimensional subspace $\ag_1 \subset \ag$ such
that $\g_1 = \ag_1 \oplus \n$, with the induced inner product, is again a standard Einstein metric solvable Lie algebra.
What is more, $\g_1$ is essentially uniquely determined by $\n$, and all the standard Einstein metric solvable Lie
algebras with the nilradical $\n$ can be obtained from $\g_1$ via a known procedure (by adjoining appropriate
derivations).

In particular, the geometry (and the algebra) of a standard Einstein metric solvable Lie algebra is completely encoded
in its nilradical (which coincides with the derived algebra). A nilpotent Lie algebra which can be a nilradical of a
standard Einstein metric solvable Lie algebra is called an \emph{Einstein nilradical}. The up-to-date list of Einstein
nilradicals can be found in \cite{L4}.

In this paper we study Einstein metric solvable Lie algebras whose nilradical is the largest possible nilpotent Lie
algebra of the given nilpotency class, a free Lie algebra. Our motivation here is twofold. Firstly, the majority
of examples (and counterexamples) of Einstein nilradicals are either two-step nilpotent or filiform (see, however
\cite{P} and \cite{T}).
Secondly, despite the fact that the algebraic and the combinatorial properties of free Lie algebras
were extensively studied for a long time by many authors (see \cite{B, R}), there seem to be not many results on their
geometry in the literature.

As a byproduct, we show that
the Einstein condition imposes the restrictions on the
growth of the dimension
of the terms in the descending central series. As a result, only six free Lie algebra of nilpotency class higher
than two are Einstein nilradicals (see the Theorem below).

For $m \ge 2, \; p \ge 1$, denote $\f(m, p)$ the (unique up to an isomorphism) free Lie algebra, which is $p$-step and
has $m$ generators. A free Lie algebra $\f(m, p)$ admits a canonical $\mathbb{N}$-gradation:
$\f(m, p) = \oplus_{k=1}^p \p(m, k)$, where $\p(m, k)$ is the space of \emph{Lie polynomials} of degree $k$ on $m$
generators. The $\p(m, k)$'s are the eigenspaces of the \emph{canonical derivation} $\Phi$ of $\f(m, p)$
defined by $\Phi_{|\p(m, k)} = k \, \id_{|\p(m, k)}$.
With a basis $e_1, \ldots, e_m$ of generators (that is, a basis for $\p(m,1)$) fixed,
we denote $e_{i_1i_2\ldots i_k}=[[\ldots[e_{i_1}, e_{i_2}],\ldots],e_{i_k}] \in \p(m, k), \; k \le p$. The space
$\p(m,k)$ has a basis of vector of this form.
For a vector $e_{i_1i_2\ldots i_k}$, define its \emph{content} to be
an $m$-dimensional vector $c(e_{i_1i_2\ldots i_k})$ whose $i$-th component is the number of occurrences of the label
$i$ among $i_1, \ldots, i_k$.

Introduce an inner product on $\f(m, p)$ in such a way that the spaces $\p(m,k)$ are orthogonal, $\<e_i, e_j\> = \K_{ij}$
for some fixed basis $\{e_i\}$ for $\p(m,1)$, and $\<e_{ij}, e_{ls}\> = \K_{il}\K_{js}$, for $i < j, \, l < s$.
We also require that any two vectors $e_{i_1i_2\ldots i_k}$ different by content are orthogonal and that the
group of permutations of $\{1, 2, \ldots, m\}$ in the subscripts of $e_{i_1i_2\ldots i_k}$ acts by isometries.
Any inner product on $\f(m, p)$ satisfying these conditions is called \emph{admissible}. An admissible
inner product on $\f(m, p)$ is unique (up to an automorphism), when $p = 1, 2$, but there are many nonisometric
admissible inner products for $p > 2$.

We prove the following theorem:

\begin{theorem*}\label{th:f}
Let $\g$ be a solvable Lie algebra with the nilradical $\n= \f(m, p)$.

1. The algebra $\g$ admits an Einstein inner product only if

(a) $p = 1$, that is, the nilradical is abelian.

(b) $p = 2$.

(c) $p=3, \; m=2,3,4,5$.

(d) $p=4, 5, \; m =2$.

In all the cases, the Einstein inner product on $\g$ is standard and $\dim \g \le \dim \n + m$.

2. Suppose additionally that $\g$ is a one-dimensional extension of its nilradical,
$\g = \mathbb{R} H \oplus \n$. Then $\ad_{H|\n} = \hat c \, \Phi$, where $\Phi$ is the canonical derivation,
$\hat c = -c \, \Tr \Phi \, (\Tr \Phi^2)^{-1}$, and $c \dim \g$ is the scalar curvature.
Up to an isometry and scaling, the Einstein inner product on $\g$ is defined as follows:
$H \perp \n$, $\|H\|^2 = \hat c \, \Tr \, \Phi$, the inner product on $\f(m, p)$ is admissible.

In case (c), define the inner product on $\p(m, 3)$ by
$$
\|e_{iji}\|^2 = 3 t, \; \|e_{ijk}\|^2 = 2 t, \; \<e_{ijk}, e_{jki}\> = -t, \quad
\text{where} \; t = \frac{m+1}{-m^2+4m+8},
$$
for $i, j, k$ pairwise distinct. % xi^2 = 3t

In case (d), for $\f(2, 4)$, define the inner product on $\p(2, 3)$ and $\p(2, 4)$ by
\begin{equation*}
\|e_{121}\|^2 = \tfrac94 \, , \quad \|e_{1211}\|^2 = \tfrac92 \, , \quad \|e_{1212}\|^2 = \tfrac94 \, .
\end{equation*}

In case (d), for $\f(2, 5)$, define the inner product on $\p(2, 3), \; \p(2, 4)$ and $\p(2, 5)$ by
\begin{gather*}
\|e_{121}\|^2 = \xi^2, \quad
\|e_{1211}\|^2 = \sigma^2 , \quad \|e_{1212}\|^2 = \sigma^2/2, \\
\|e_{12111}\|^2 = \alpha^2, \quad \|u\|^2 = 3 \alpha^2, \quad
\|[e_{121}, e_{12}]\|^2 = \gamma^2, \quad \<u, [e_{121}, e_{12}]\> = 0,
\end{gather*}
where $u = e_{12211} + e_{12121} + e_{12112}$ and
$\xi^2 = 54 t, \; \sigma^2 = 3/4+375 t, \; \alpha^2 = 76 t \sigma^2, \; \gamma^2 = 27 t (1 + 128 t)$,
$t=(131+5\sqrt{745})/2928$.
\end{theorem*}

\begin{remark}
The fact that abelian and two-step free Lie algebras are Einstein nilradicals is well known. Standard Einstein metric
solvable Lie algebras with an abelian nilradical were classified in \cite{A2}. Any two-step free Lie algebra is an
Einstein nilradical, as follows from the result of \cite{GK}.
\end{remark}

\begin{remark}
The fact that a metric solvable Einstein Lie algebra with free nilradical is always standard is proved in
\cite[Theorem 5]{Ni}. This is why our main focus throughout the paper is on the Einstein nilradicals,
the nilradicals of standard Einstein metric solvable Lie algebras.
\end{remark}

The paper is organized as follows. In Section~\ref{s:facts} we provide the necessary background on Einstein solvmanifolds
(\ref{ss:einstein}) and on free Lie algebras (\ref{ss:free}). The proof of the Theorem (Section~\ref{s:proof}) is split
into three parts. In \ref{ss:convex} we obtain some general inequality-type restrictions on the eigenvalue
type of an Einstein nilradical. Using that we show in \ref{ss:der} that only a small number of the algebras
$\f(m, p)$ with $p > 2$ can be Einstein nilradicals. In \ref{ss:metric} we compute the inner product and the Ricci
tensor for them, hence completing the proof.

The author is grateful to W.Ziller and to J.Lauret for useful comments and references.

\section{Preliminaries}
\label{s:facts}

\subsection{Einstein solvmanifolds}
\label{ss:einstein}

For an inner product $\< \cdot, \cdot \>$ on a Lie algebra $\g$, define the \emph{mean curvature vector} $H$ by
$\<H, X\> = \Tr \ad_X$ (clearly, $H$ is orthogonal to the derived algebra of $\g$). For $A \in \End(\g)$, let $A^*$
be its metric adjoint and $S(A) =\frac12 (A +A^*)$ be the symmetric part of $A$. Let $\Ric$ be the Ricci
$(0, 2)$-tensor (a quadratic form) of $(\g, \< \cdot, \cdot \>)$, and $\ric$ be the \emph{Ricci operator}, the symmetric
operator associated to $\Ric$.

The Ricci operator of $(\g, \< \cdot, \cdot \>)$ is implicitly defined by
\begin{equation}\label{eq:riccidef}
\Tr \Bigl(\ric + S(\ad_H) + \frac12 B \Bigr)
\circ A = \frac14 \sum_{i,j} \<A[E_i, E_j] - [AE_i, E_j] - [E_i, AE_j], [E_i, E_j]\>,
\end{equation}
for any $A \in \End(\g)$, where $\{E_i\}$ is an orthonormal basis for $\g$, and $B$ is the symmetric operator associated
to the Killing form of $\g$. Note that the right-hand side contains the action of the coboundary operator on $A$.

If $(\n, \< \cdot, \cdot \>)$ is a nilpotent metric Lie algebra, then $H = 0$ and $B = 0$, so \eqref{eq:riccidef} gives
\begin{equation}\label{eq:riccinil}
\Tr (\ric_{\n} \circ A) = \frac14 \sum_{i,j} \<A[E_i, E_j] - [AE_i, E_j] - [E_i, AE_j], [E_i, E_j]\>,
\end{equation}
or explicitly, for a vector $X \in \n$,
\begin{equation}\label{eq:Riccinil}
\Ric_{\n} (X) = \frac14 \sum_{i,j} \<X, [E_i, E_j]\>^2 - \frac12 \sum_{i} \| [X, E_i]\|^2.
\end{equation}

\begin{definition*} \cite{H}
An inner product on a solvable Lie algebra $\g$ is called \emph{standard}, if the orthogonal complement to
the derived algebra $[\g, \g]$ is abelian. A metric solvable Lie algebra $(\g, \<\cdot,\cdot\>)$ is called
\emph{standard}, if the inner product $\<\cdot,\cdot\>$ is standard.
\end{definition*}

As it is proved in \cite{AK}, any Ricci-flat metric solvable Lie algebra is flat. By the result of \cite{DM},
any Einstein metric solvable unimodular Lie algebra is also flat. In what follows, we always assume $\g$
to be nonunimodular ($H \ne 0$), with an inner product of a strictly negative scalar curvature $c \dim \g$.

A standard metric Einstein Lie algebra admits a rank-one reduction \cite[Theorem 4.18]{H}. This means that if
$(\g, \< \cdot, \cdot\>)$ is such an algebra, with the nilradical $\n$ and the mean curvature vector $H$, then the
subalgebra $\g_1 = \mathbb{R}H \oplus \n$, with the induced inner product, is also Einstein and standard. What is
more, the derivation $\phi=\ad_{H|\n}:\n \to \n$ is symmetric with respect to the inner product, and all its
eigenvalues belong to $\a \mathbb{N}$ for some constant $\a > 0$. This implies, in particular, that the nilradical $\n$
of a standard Einstein metric solvable Lie algebra admits an $\mathbb{N}$-gradation defined by the eigenspaces of $\phi$.
As it is proved in \cite[Theorem 3.7]{L1}, a necessary and sufficient condition for a metric nilpotent algebra
$(\n, \< \cdot, \cdot\>)$ to be the nilradical of a standard metric Einstein solvable Lie algebra is
\begin{equation}\label{eq:ricn}
    \ric_\n = c \, \id_\n + \phi,
\end{equation}
where $c \dim \g < 0$ is the scalar curvature of $(\g, \< \cdot, \cdot\>)$. This equation, in fact, defines
$(\g, \< \cdot, \cdot\>)$ in the following sense: given a metric nilpotent Lie algebra whose Ricci operator
satisfies \eqref{eq:ricn}, with some constant $c < 0$ and some $\phi \in \Der(\n)$, one can define $\g$ as a
one-dimensional extension of $\n$ by $\phi$. For such an extension $\g = \mathbb{R}H \oplus \n, \; \ad_{H|\n} = \phi$,
and the inner product
defined by $\<H, \n \> = 0,\; \|H\|^2 = \Tr \phi$ (and coinciding with the existing one on $\n$) is Einstein, with the
scalar curvature $c \dim \g$. Following \cite{L1} we call a nilpotent Lie algebra $\n$ which admits an inner product
$\< \cdot, \cdot\>$ and a derivation $\phi$ satisfying \eqref{eq:ricn} an \emph{Einstein nilradical}, the
corresponding derivation $\phi$ is called an \emph{Einstein derivation}, and the inner product $\< \cdot, \cdot\>$
the \emph{nilsoliton metric}.

As it is proved in \cite[Theorem 3.5]{L1}, a nilpotent Lie algebra admits no more than one nilsoliton metric, up to
conjugation and scaling (and hence, an Einstein derivation, if it exists, is unique, up to conjugation and scaling).

Equation \eqref{eq:ricn}, together with \eqref{eq:riccinil}, implies that
\begin{equation}\label{eq:trace}
    \Tr (\phi \circ \psi) = - c \, \Tr \psi,  \quad \text{for any $\psi \in \Der(\n)$}.
\end{equation}
Substituting $\psi = \phi$ in this equation, one finds that $-c = \Tr(\phi^2)/\Tr(\phi)$. If
$\mu_1 < \ldots < \mu_p$ are the eigenvalues of $\phi$, with $d_1, \ldots, d_p$ the corresponding multiplicities, we
obtain $ -c \sum \mu_k d_k = \sum \mu_k d_k^2$. The set $(\mu_1 < \ldots < \mu_p; \; d_1, \ldots, d_p)$ is called
the \emph{eigenvalue type} of the Einstein metric solvable Lie algebra $(\g, \< \cdot, \cdot\>)$. With some abuse of
language, we will call the same set $(\mu_1 < \ldots < \mu_p; \; d_1, \ldots, d_p)$ the
\emph{eigenvalue type} of $\phi$.

\subsection{Free Lie algebras}
\label{ss:free}

In this section, we collect some known facts on the free Lie algebras (for details see \cite{B, R}).

From among $p$-step nilpotent Lie algebras on $m \ge 2$ generators, a \emph{free Lie algebra} $\f(m, p)$ is the one
having the maximal dimension ($\f(m, p)$ is unique, up to an isomorphism). In other words, given the generators
$e_1, \ldots, e_m$, the algebra $\f(m, p)$ is the linear span of all the $k$-folded brackets of the $e_i$'s, $k \le p$,
with the only relations between these brackets coming from the skew-symmetricity and the Jacobi identity. For every
$k=1, \ldots, p$, the subspace of $\f(m, p)$ spanned by the $k$-folded brackets is the space of \emph{Lie polynomials}
$\p(m, k)$. In particular, $\p(m, 1) = \Span (e_1, \ldots, e_m)$. The direct sum decomposition
$\f(m, p) = \oplus_{k=1}^p \p(m, k)$ is an $\mathbb{N}$-gradation. It corresponds to the \emph{canonical derivation}
$\Phi$ defined by
\begin{equation}\label{eq:canPhi}
    \Phi \Bigl(\sum\nolimits_{k=1}^p X_k\Bigr) = \sum\nolimits_{k=1}^p k X_k, \quad \text{where } X_k \in \p(m, k).
\end{equation}

The dimension of the space $\p(m, k)$ is given by
\begin{equation}\label{eq:dimfree}
    d_k(m) = \dim \p(m, k) = k^{-1} \sum\nolimits_{d|k} \mu(d) m^{k/d},
\end{equation}
where the sum is taken over all the factors $d$ of $k$, and $\mu(d)$, the M\"{o}bius function, is defined by
$\mu(n) = 0$, if $n$ is divisible by a square of a prime number, and $\mu(n) = (-1)^{\omega(n)}$ otherwise, where
$\omega(n)$ is the number of prime factors of $n$ (in particular, $\mu(1) = 1$). The table below lists
the dimensions $d_k(m)$ for $k \le 7, \; m \le 3$, and for $k \le 14, \; m = 2$:

\begin{table}[h]
\begin{center}
\begin{tabular}{|c|c|c|c|}
  \hline
  $k$ & $d_k(m)$ & $d_k(2)$ & $d_k(3)$ \\
  \hline
  $1$ & $m$ & 2 & 3 \\
  $2$ & $(m^2-m)/2$ & 1 & 3 \\
  $3$ & $(m^3-m)/3$ & 2 & 8 \\
  $4$ & $(m^4-m^2)/4$ & 3 & 18 \\
  $5$ & $(m^5-m)/5$ & 6 & 48 \\
  $6$ & $(m^6 - m^3 - m^2 + m)/6$ & 9 & 116 \\
  $7$ & $(m^7-m)/7$ & 18 & 312 \\
  \hline
\end{tabular}
\hskip 1.5cm
\begin{tabular}{|c|c|}
  \hline
  $k$ & $d_k(2)$ \\
  \hline
  $8$ & $30$ \\
  $9$ & $56$ \\
  $10$ & $99$\\
  $11$ & $186$ \\
  $12$ & $335$ \\
  $13$ & $630$ \\
  $14$ & $1161$\\
  \hline
\end{tabular}
\caption{Dimensions of the spaces of Lie polynomials.}\label{table}
\end{center}
\end{table}

With the generators $e_1, \ldots, e_m$ fixed, we denote
$e_{i_1i_2\ldots i_k}=[[\ldots[e_{i_1}, e_{i_2}],\ldots],e_{i_k}] \in \p(m, k)$, for $k \le p$. Every space
$\p(m,k)$ is spanned by the vectors $e_{i_1i_2\ldots i_k}$, but they are not in general linearly independent (e.g.,
$e_{1221} = e_{1212}$). There are several convenient bases for $\f(m,p)$, one of which is the Hall
basis \cite[\S 2.11]{B}. As the direct computations involving a particular basis will be needed only in few
low-dimensional cases, we choose a more symmetric basis, partially consisting of the $e_{i_1i_2\ldots i_k}$'s
and partially of the elements of the Hall basis (see Table~\ref{table2} in Lemma~\ref{l:product}).
For a vector $e_{i_1i_2\ldots i_k}$, we define its \emph{content} to be
an $m$-dimensional vector $c(e_{i_1i_2\ldots i_k})$, whose $i$-th component is the number of occurrences of the label
$i$ among $i_1, \ldots, i_k$.

For a free Lie algebra, any assignment of the images to the generators extends to a (unique) derivation: for any linear
map $L: \p(m,1) \to \f(m,p)$, there exists a unique derivation whose restriction to $\p(m,1)$ coincides with $L$. In
particular, any endomorphism $L$ of $\p(m,1)$ extends to a derivation $\rho(L) \in \Der (\f(m,p))$. The spaces $\p(m,k)$
are
invariant with respect to $\rho(L)$. For every $k = 1, \ldots , p$, let $\rho_k(L)$ be the restriction of $\rho(L)$ to
the $\p(m,k)$. Then $\rho_k$ is a representation of the Lie algebra $\mathfrak{gl}(m)$ on the space $\p(m,k)$ of Lie
polynomials. Similarly, the action of any $S \in \mathrm{GL}(m)$ on $\p(m,1)$ extends to an automorphism $R(S)$
of $\f(m,p)$, with the spaces $\p(m,k)$ invariant with respect to $R(S)$. The restriction $R_k(S)$ of
$R(S)$ to $\p(m,k)$ defines a representation of $\mathrm{GL}(m)$ on the space $\p(m,k)$ (with
$dR_k = \rho_k$). Note that
$R_k$ is, in general, reducible. For example, the representation $R_6$ of $\mathrm{GL}(2)$ on $\p(2, 6)$ has a
one-dimensional invariant subspace spanned by $I = [e_{121}, e_{122}]$ (such a vector is called a \emph{Lie invariant}).

\section{Proof of the Theorem}
\label{s:proof}

By \cite[Theorem 5]{Ni}, any metric solvable Einstein Lie algebra with free nilradical is standard. Therefore to
prove the Theorem we have to determine which free Lie algebras are Einstein nilradicals and then to find their
one-dimensional Einstein metric solvable extensions.

One of the most efficient methods of finding Einstein nilradicals
is the variational one based on the study of the Ricci flow on the space of Lie brackets on the given
linear space $\Rn$ with a fixed inner product (see \cite{L1, L2, L3, LW}). The limit of the flow (which always exists and
is a critical point of the normalized squared norm of the Ricci tensor) is a Lie bracket on the given inner product
space, which turns it into an Einstein
nilradical. If the limiting bracket $\lambda$ lies in the $\GL(n)$-orbit of the initial one, $\mu$, then the nilpotent
Lie algebra $(\Rn, \mu)$ is an Einstein nilradical. If not (the bracket $\mu$ degenerates to a nonisomorphic
bracket $\lambda$), one can prove that $(\Rn, \mu)$ is not an Einstein nilradical using the stratification on
the space of skew-symmetric bilinear maps from $\Rn$ to itself given in \cite{LW}.

Another method, introduced in \cite{P} gives an elegant and very computable criterion of finding Einstein nilradicals
using the Gramm matrix of the root vectors with respect to the basis of Ricci eigenvectors.

These methods proved to be extremely successful both in constructing
Einstein nilradicals and in showing that a given nilpotent Lie algebra is not an Einstein nilradical.

It seems, however, that one can hardly find a basis for free Lie algebra, which is sufficiently well adapted for
computing the Ricci tensor, or whose structural constants have a nice form.
Even with a very neatly chosen basis (Section~\ref{ss:metric}), the explicit
calculation of the Ricci tensor is rather unpleasant.

We use a somewhat different approach starting with the Lie algebra itself and then constructing the inner product.
It consists of the following three steps.

First, we show that the eigenvalue type of an (arbitrary) Einstein solvmanifold satisfy a set of inequality-type
conditions (Lemma~\ref{l:convex}). That lemma is independent of the assumption of the freeness of the nilradical and
may represent a certain interest by itself
(as an application of Lemma~\ref{l:convex}, we consider possible eigenvalue types of Einstein
nilradicals whose Einstein derivation has three eigenvalues). One of the consequences of Lemma~\ref{l:convex} is that
the dimension of the eigenspaces of the Einstein derivation cannot grow too fast. This is the key idea of the proof,
the rest is in the sense technical.

Secondly, using \eqref{eq:trace} we find the Einstein derivation (if it exists) of a free Lie algebra. As expected,
the Einstein derivation is a multiple of the canonical derivation $\Phi$, up to conjugation
(Lemma~\ref{l:der}). As the dimension of the eigenspaces of $\Phi$ grows fast enough, namely $d_k(m) \sim m^k/k$
by equation \eqref{eq:dimfree}, the inequalities of Lemma~\ref{l:convex} imply that only finitely many free Lie
algebras, other than the abelian and the two-step ones, can be Einstein nilradicals (they are listed in
Lemma~\ref{l:nonexist}).

On the third step, for each of the free algebras not eliminated by Lemma~\ref{l:nonexist}, we explicitly construct the
unique inner product which has a
chance to be a nilsoliton one (Lemma~\ref{l:product}). The construction is based on the following: firstly, the
Einstein derivation must be symmetric, and secondly, for any derivation $\psi$ which commutes with the Einstein
derivation,
its symmetric part must also be a derivation (this follows from \cite[Theorem 4.12]{H} and the fact that any such $\psi$
extends to a derivation of the rank-one Einstein solvable extension of $\n$).
To finish the proof, we then find the Ricci tensor in the remaining cases and check, whether it
satisfies \eqref{eq:Riccinil}.

\subsection{Eigenvalue type}
\label{ss:convex}

Let a nilpotent Lie algebra $\n$ be an Einstein nilradical, with $\phi$ the corresponding Einstein derivation whose
eigenvalue type is $(\mu_1 < \ldots < \mu_p; \; d_1, \ldots, d_p)$. Recall that all the $\mu_i$'s are positive and
$\mu_i/\mu_j \in \mathbb{Q}$. Let $\n_1, \ldots, \n_p$ be the eigenspaces of $\phi$ corresponding to the eigenvalues
$\mu_1, \ldots, \mu_p$, respectively, $\dim \n_i = d_i$. In \cite[Theorem 4.14]{H} it is shown that the numbers
$\mu_i, d_i$ satisfy a certain system of equations arising from the following. In the Euclidean space $\mathbb{R}^p$
with an orthonormal basis
$f_1, \ldots, f_p$, define a subset $F = \{f_k-f_i-f_j \,| \, \mu_i + \mu_j = \mu_k \}$. Then any vector
$\nu =(\nu_1, \ldots, \nu_p)^t \in \mathbb{R}^p$, which is orthogonal to $\Span(F)$, defines a derivation of $\n$
acting by multiplication by $\nu_i$ on every subspace $\n_i$. As for any such derivation $\psi$, equation
\eqref{eq:trace} must be satisfied, one gets a system of equations for the $\mu_i$'s and the $d_i$'s.

We start with the following lemma which adds some further inequality-type restrictions on the $d_i$'s. They
seem to be most effective in the cases when no or few equations are present (that is, when $\mathrm{codim}\,\Span(F)$
is small, for instance, when $\mu_i = i$).

\begin{lemma} \label{l:convex}
Let $\n$ be an Einstein nilradical, with $\phi$ the corresponding Einstein derivation
whose eigenvalue type is $(\mu_1 < \ldots < \mu_p; \; d_1, \ldots, d_p)$.
In a Euclidean space $\mathbb{R}^p$ with a fixed orthonormal basis
$f_1, \ldots, f_p$, define a subset $F = \{f_k-f_i-f_j \,| \, \mu_i + \mu_j = \mu_k \}$ and the vectors
$[\mu] = (\mu_1, \ldots, \mu_p)^t$, and $[1] = (1, \ldots, 1)^t$. Let
$D$ be a diagonal $p \times p$-matrix, whose diagonal entries are $d_1, \ldots, d_p$, respectively.

Then the vector $v = (v_1, \ldots, v_p)^t \in \mathbb{R}^p$ with the components
\begin{equation}\label{eq:v}
v=D(\Tr(\phi) [\mu]- \Tr(\phi^2) [1])
\end{equation}
lies in the convex cone spanned by $F$.
\end{lemma}

Explicitly, the components of $v$ are given by $v_k = d_k \sum_i d_i \mu_i(\mu_k-\mu_i)$.
Note that both the vector $v$ and all the vectors from $F$ lie in the orthogonal complement to $[\mu]$.

The idea of the convex hull already appeared in the study of Einstein nilradicals \cite[Section~3]{LW},
here we apply it to the eigenvalue type.

\begin{proof} Let $\n_i, \; i=1, \ldots, p$, be the eigenspace of $\phi$ corresponding to the eigenvalue $\mu_i$.
For $1 \le i, j \le p$, define the numbers $T_{ij}=\|[\n_i, \n_j]\|^2 = \sum_{\a \beta} \|[E_\a, E_\beta]\|^2$,
where $\{E_\a\}$ and $\{E_\beta\}$ are orthonormal bases for $\n_i$ and $\n_j$ respectively. Clearly, $T_{ij} \ge 0$
and $T_{ij} = T_{ji}$.

For every $k=1, \ldots, p$, let $A_k$ be an endomorphism of $\n$ whose restriction to $\n_k$ is the identity, and
which sends to zero all the eigenspaces $\n_i$ with $i\ne k$. Substituting $A_k$ as $A$ into \eqref{eq:riccinil}
and using \eqref{eq:ricn}, we get
$$
\sum_{i,j: \mu_i + \mu_j = \mu_k} T_{ij} - 2 \sum_j T_{kj} = 4 (\Tr(\phi))^{-1} v_k, \qquad k = 1, \ldots, p.
$$
This is a system of linear equations for the nonnegative numbers $T_{ij}, \, i \le j$.
The nonzero columns of its matrix are precisely the vectors from the set $F$ (some of them doubled). This proves
the lemma, as all the $T_{ij}$'s are nonnegative and $\Tr(\phi) >0$.
\end{proof}

Before applying Lemma~\ref{l:convex} to the free algebras (see Lemma~\ref{l:nonexist} below), we
illustrate the sort of restrictions imposed by it on the eigenvalue type of an Einstein derivation with a small number
of eigenspaces (note that Lemma~\ref{l:convex} says nothing if $p = 2$).

\begin{example}
Consider the case $p=3$. Let $\mu_1 < \mu_2 < \mu_3$ be the eigenvalues
of $\phi$, with the multiplicities $d_1, d_2, d_3$, respectively.
We will use the result of \cite[Proposition 3.3]{L3}, from which it follows that an abelian factor of $\n$ (if any)
is an eigenspace of $\phi$ whose eigenvalue cannot be the maximal or the minimal one among the $\mu_i$'s.

The set $F \subset \mathbb{R}^3$ contains at most two vectors. Clearly, $F$ cannot be empty (otherwise $\n$ is abelian
and $\phi$ does not satisfy \eqref{eq:trace}).
If $\#(F) = 1$, two cases are possible:
either $\mu_1 + \mu_2 = \mu_3$ (and $\mu_2 \ne 2 \mu_1$), or $2\mu_1 = \mu_3$. The first case, by 5.3(D) of
\cite{H}, gives $\mu_2 d_1 (d_2 + d_3) = \mu_1 d_2 (d_1 + d_3)$, Lemma~\ref{l:convex} imposes no further restrictions.
In the second one, $\n_2$ is an abelian factor and one has $\mu_2 (d_1 + 2d_3) = \mu_1 (d_1 + 4d_3)$
(see \cite{L3}). The corresponding rank-one Einstein solvmanifold (if it exists) is the  rank-one reduction of the
Riemannian product of a two-step Einstein solvmanifold with the eigenvalue type $(1<2; \, d_1, d_3)$ and the hyperbolic
space $H^{d_2}$ of the same Ricci curvature.

Assume now that $\#(F)=2$. Then $(\mu_1, \mu_2, \mu_3)$ is proportional either to $(1, 2, 3)$ or to $(1, 2, 4)$. In the
first case, Lemma~\ref{l:convex} gives
\begin{equation}\label{eq:convex3}
    d_3^{-1} + 2 d_2^{-1} \ge d_1^{-1}.
\end{equation}
The second one is not possible. Indeed, let $[\n_1, \n_1] = \n_2^1 \subset \n_2$ and let $\n_2^2$ be a subspace of $\n_2$
complementary to $\n_2^1$. Then $\n_2^1$ lies in the center of $\n$ (from the Jacobi identity) and is nonzero (otherwise
$\n_1$ is an abelian factor, which is impossible, as $1$ is the smallest eigenvalue). Similarly, $\n_2^2$ is nonzero, as
otherwise $\n_3$ would be an abelian factor. Let $d_{21} = \dim \n_2^1,\; d_{22} = \dim \n_2^2$.
For any real numbers $a, b$, the endomorphism $\psi$, which acts as a
multiplication by $a, 2a, b, 2b$ on the subspaces $\n_1, \n_2^1, \n_2^2, \n_3$, respectively, is a derivation.
Substituting
such a $\psi$ into \eqref{eq:trace} we get $(d_1 + 4d_{21})/(d_1 + 2d_{21}) = (2d_{22} + 8d_3)/(d_{22} + 2d_3)$, which
is a contradiction, as the left-hand side is less than $2$, but the right-hand side is bigger than $2$.
\end{example}
\begin{example} For the eigenvalue type $(1 < 2 < 3 < 4; \, d_1, d_2, d_3, d_4)$, Lemma~\ref{l:convex} gives the
following inequalities:
\begin{equation}\label{eq:convex4}
\begin{split}
6 d_4 d_1 + 2 d_4 d_3 + d_2 d_1  + 2 d_3 d_1 &\ge d_2 d_3,\\
6 d_4 d_1 + d_2 d_1 + d_2 d_3 + 4 d_3 d_1 &\ge 2 d_4 d_3,\\
3 d_4 d_1 + 4 d_4 d_2 + 2 d_3 d_1 + 2 d_2 d_3 &\ge d_4 d_3.
\end{split}
\end{equation}

\end{example}

From Lemma~\ref{l:convex} it follows, in particular, that the dimensions $d_k$ cannot grow very fast. This is the
reason why only finitely many free algebras, apart from the abelian and the two-step ones, can be Einstein
nilradicals (as according to \eqref{eq:dimfree}, $d_k(m)\sim m^k/k$).

\subsection{Einstein Derivation}
\label{ss:der}

Although a free Lie algebra has many semisimple derivations whose eigenvalues are natural numbers (hence many
nonisomorphic $\mathbb{N}$-gradations), only one of them, up to conjugation and scaling, can be an Einstein
derivation $\phi$.
The uniqueness immediately follows from the result of \cite[Theorem 3.5]{L1} on the uniqueness of the nilsoliton
metric. The following lemma shows that $\phi$ is proportional to the obvious candidate: the canonical derivation $\Phi$
defined by \eqref{eq:canPhi}.

\begin{lemma}\label{l:der}
An Einstein derivation $\phi$ of $\n =\f(m, p)$, up to conjugation by an automorphism of $\f(m, p)$, is given by
$$
\phi = \hat c \, \Phi,
$$
where $\hat c = -c \, \Tr \Phi \, (\Tr \Phi^2)^{-1}$, and $c (\dim \n +1)$ is the scalar curvature of the
Einstein inner product on the rank-one extension of $\f(m, p)$ by $\phi$.
\end{lemma}
\begin{proof}
The uniqueness of $\phi$, up to conjugation, follows from the uniqueness of the nilsoliton metric. If an Einstein
derivation $\phi$ exists, it has to satisfy \eqref{eq:trace}. By \cite[Proposition~1]{Ni}, a semisimple derivation
with real eigenvalues satisfying \eqref{eq:trace} is unique up to conjugation and scaling,
so to prove the lemma it suffices to show that $\hat c \, \Phi$
satisfies \eqref{eq:trace}. As $\Phi$ is semisimple, the derivation $\ad_\Phi$ of $\Der(\n)$ is also semisimple.
If $\psi \in \Der(\n)$ is an eigenvector of $\ad_\Phi$ with a nonzero eigenvalue, then
$\Tr \, \Phi \circ \psi = \Tr \, \psi = 0$, and \eqref{eq:trace} is obviously satisfied. So it is
sufficient to consider only those $\psi$ which commute with $\Phi$. For any such $\psi$, the spaces $\p(m, k)$ are
invariant, and moreover, $\psi = \rho(L)$ for some endomorphism $L$ of $\p(m,1)$. As $\rho(\id) = \Phi$, it
suffices to prove the following: for any $L$ with $\Tr L = 0$ and for any $k \le p, \quad \Tr \rho_k(L) = 0$. This
follows from the fact that $\rho_{k|\mathfrak{sl}(m)}$ is a representation of the simple algebra $\mathfrak{sl}(m)$.
\end{proof}

By Lemma~\ref{l:der}, if a free Lie algebra $\f(m, p)$ is an Einstein nilradical, then the eigenvalue type
of its Einstein derivation is $(\hat c < 2 \, \hat c < \ldots < p \, \hat c ; \; d_1(m), d_2(m), \ldots, d_p(m))$, where
the numbers $d_k(m)$ are given by \eqref{eq:dimfree}.

\begin{lemma}\label{l:nonexist}
A free Lie algebra $\f(m, p)$ with $m \ge 2,\; p \ge 3$ can be an Einstein nilradical only in the following cases:

\emph{(a)} $m = 2, \; p \le 7$;

\emph{(b)} $(m, p) = (3, 3), \, (4, 3), \, (5, 3),\, (3, 4)$.
\end{lemma}

As it will be shown in Section~\ref{ss:metric}, not all the free algebras listed above are Einstein
nilradicals.

\begin{proof}
We start by considering the cases $p=3$ and $p=4$. By \eqref{eq:convex3}, for $\f(m, 3)$ to be an Einstein nilradical,
the equation $d_3(m)^{-1} + 2 d_2(m)^{-1} \ge d_1(m)^{-1}$ must hold. Substituting $d_1(m), \, d_2(m), \, d_3(m)$ from
Table~\ref{table} we get $-m^2 + 4m + 8 \ge 0$, so $m \le 5$. When $p = 4$, inequalities \eqref{eq:convex4} have to
be satisfied. Substituting $d_1(m), \, d_2(m), \, d_3(m), \, d_4(m)$ into the second one we get
$m^4 + m^3-11 m^2- 18 m - 10 < 0$, which is only true when $m = 2, 3$.

From now on, we assume that $p \ge 5$.

By Lemma~\ref{l:convex}, to prove that $\f(m, p)$ is not an Einstein nilradical, it suffices to produce
a vector $a \in \mathbb{R}^p$ such that $a \perp v$ and all the scalar products of $a$ with the vectors from $F$ are
negative. For $\f(m, p)$, the components of the vector $v$ are given by
$v_k = d_k(m) ( k \sum_{j=1}^p j d_j(m) - \sum_{j=1}^p j^2 d_j(m))$ (equation \eqref{eq:v}). The set $F$ consists of
all the vectors of the form $f_{i+j}-f_i-f_j$, with $i+j \le p$.

In the most cases, such a vector $a$ can be taken as
\begin{equation*}
    a = \Bigr(\frac{1}{d_1(m)}, \frac{1}{d_2(m)}, \ldots, \frac{1}{d_{p-2}(m)}, -\frac{(p-2)(p+1)}{2 d_{p-1}(m)},
\frac{(p-2)(p-1)}{2 d_p(m)}\Bigl)^t.
\end{equation*}
Then $a \perp v$, and to check that the scalar product of $a$ with all the vectors from $F$ is negative, we have to show
that $a_{i+j} < a_i + a_j$, for all $i, j$ with $i+j \le p$, which gives the following inequalities:
\begin{align}
&\frac{1}{d_i(m)} + \frac{1}{d_j(m)} > \frac{1}{d_{i+j}(m)} \, , & \quad i + j \le p-2, \label{eq:d1} \\
&\frac{1}{d_i(m)} + \frac{1}{d_j(m)} > \frac{(p-2)(p-1)}{2 d_p(m)} \, , & \quad i + j = p, \; i, j < p-1, \label{eq:d2}\\
&\frac{1}{d_1(m)} - \frac{(p-2)(p+1)}{2 d_{p-1}(m)} > \frac{(p-2)(p-1)}{2 d_p(m)} \label{eq:d3} \, .
\end{align}

Note that inequality \eqref{eq:d1} becomes an equality when $i=j=1, \, m=2$, so $a \perp f_2 - 2f_1$. However, if the
scalar product of $a$ with all the other vectors from $F$ is negative, the only possible way for $v$ to lie in the convex
cone spanned by $F$ is to be (positively) proportional to $f_2 - 2f_1$, which is clearly not the case when $p \ge 5$.

We make use of the following three properties of the $d_k(m)$'s:
\begin{enumerate}[(i)]
    \item $d_k(m)$ increases as a function of $m$; $d_k(m)$ increases as a function of $k$, with the only exception:
    $d_1(2) > d_2(2)$.
    \item $d_k(m) \le k^{-1} m^k$, for all $m \ge 2,\; k \ge 1$.
    \item $d_k(m) > k^{-1} (m^k - m^{t(k)+1}) \ge  k^{-1} (m^k - m^{[k/2]+1})$, for all $m \ge 2,\; k \ge 2$,
    where $t(k) = \max_{d \mid k,\, d < k} d$.
\end{enumerate}

The first assertion of (i) is clear: $\p(m,k)$ is a (proper) subspace of $\p(m+1,k)$. The second one follows
from the fact that $\ad_{e_1|\p(m,k)}: \p(m,k) \to \p(m,k+1)$ is injective when $k \ge 2$
and that $d_2(m) = (m^2-m)/2 \ge m = d_1(m)$ for $m > 2$ (or from (ii) and (iii) by a direct computation). Also,
from formula \eqref{eq:dimfree} we have the following estimates:
$$
m^k - m^{t(k)} - \sum\nolimits_{s=1}^{t(k)-1} m^s \le k d_k(m) \le m^k - m^{t(k)} + \sum\nolimits_{s=1}^{t(k)-1} m^s,
$$
which prove properties (ii) and (iii).

Property (i) shows that \eqref{eq:d1} is always satisfied, except for $m = 2, \; i = j = 1$, when the sides are equal.

As $d_1(m) = m$ and by property (i), to prove \eqref{eq:d3} it suffices to show that
$d_{p-1}(m) > m (p-2) p$. By property (iii), the latter inequality will be always satisfied if
$m^{p-2} - m^{[(p-1)/2]} > (p-2)(p-1)p$.
Induction arguments show that the latter inequality holds for all the pairs $(m, p)$ with $m \ge 2,\, p \ge 5$, except
for the following ones: $m = 2, \; 5 \le p \le 12$ and $(m, p) = (3, 5), (3, 6), (4, 5)$.
A direct check using the values from Table \ref{table} then shows that inequality \eqref{eq:d3} holds for all
$m \ge 2,\, p \ge 5$, unless when $(m, p) = (2,p), \, p \le 12, \; (3, 5), \; (3, 6)$.

Inequality \eqref{eq:d2} will be satisfied when $2 d_p(m) > (p-2)(p-1) d_{[p/2]}(m)$. Using the estimates from (ii)
and (iii) we find that the latter inequality is true for all $m \ge 3, \, p \ge 5$, and also for $m = 2, \; p \le 14$.
Checking directly we obtain that when $m \ge 2, \; p \ge 5$, inequality \eqref{eq:d2} only fails when $m = 2, p \le 10$.

To prove the lemma, it remains to consider the cases $(m, p) = (2,p), \, 8 \le p \le 12$, $(3, 5), \; (3, 6)$.
Take $a = (1, \ldots, 1, x)^t \in \mathbb{R}^p$, with $x$ chosen in such a way that $a \perp v$.
Then the scalar product of $a$ with every vector from $F$ is negative, provided $x < 2$. Consulting Table~\ref{table}
again we see that this is indeed the case for all the above pairs, except for $(m, p) = (2,8)$. For
that remaining case, let $a = (14, 28, 42, 5, -32, -18, -4, 10)^t$. Then $a_i + a_j \ge a_{i+j}$ for all $i,j \ge 1$ with
$i + j \le 8$, so the scalar products of $a$ with all the vectors from $F$ are nonpositive. As $\<a, v\> = 72828 > 0$,
$v$ does not lie in the convex cone spanned by $F$.
\end{proof}

\subsection{Nilsoliton metric}
\label{ss:metric}

By Lemma~\ref{l:nonexist}, a free Lie algebra $\f(m, p)$ can be an Einstein nilradical only when $p \le 2$, or when
$p = 3, \; m = 2, 3, 4, 5$, or when $p = 4,\; m =2, 3$, or when $p = 5, 6, 7,\; m =2$.

In this section, we show that all these algebras, except for $\f(3, 4), \, \f(2, 6)$ and $\f(2, 7)$, are indeed Einstein
nilradical by explicitly constructing the inner product for which \eqref{eq:Riccinil} is satisfied.

We start with the following lemma.

\begin{lemma}\label{l:productgen}
1. Let $\<\cdot,\cdot\>$ be a nilsoliton inner product on $\f(m, k)$. Then the spaces $\p(m, k)$ are orthogonal.
The inner product on $\p(m, 1)$ can be chosen arbitrarily (different choices give isometric inner products).
Choose and fix an orthonormal basis $e_1, \ldots, e_m$ of generators. Then:

\begin{enumerate}[(a)]
  \item
    The inner product on $\p(m, k)$ satisfies the equation
    \begin{equation*}
        \rho_k(L^*) = (\rho_k(L))^*,
    \end{equation*}
    for any $L \in \End(\p(m, 1))$.

  \item
  Any two vectors $e_{i_1i_2\ldots i_s}$ different by content are orthogonal.

  \item
  For any permutation $\tau \in S_m$ of the basis $e_1, \ldots, e_m$,
    \begin{equation*}
        \<e_{\tau(i_1)\ldots \tau(i_s)}, e_{\tau(j_1)\ldots \tau(j_s)}\> = \<e_{i_1\ldots i_s}, e_{j_1\ldots j_s}\>.
    \end{equation*}
\end{enumerate}

2. Let $\<\cdot,\cdot\>$ be an arbitrary inner product on $\f(m, k)$, for which the spaces $\p(m, k)$ are orthogonal and
which satisfies (a), (b) and (c). Let $\Ric$ be its $(0, 2)$ Ricci tensor. Then the spaces $\p(m, k)$ are orthogonal
with respect to the quadratic form $\Ric$, and assertions (b) and (c) remain true, if we replace $\<\cdot,\cdot\>$ by
$\Ric$.
\end{lemma}

\begin{proof}
1. The Einstein derivation $\phi$ must be symmetric, which implies that the spaces $\p(m, k)$ are orthogonal.
As the nilsoliton inner product is unique up to scaling and conjugation \cite[Theorem 3.5]{L1} and as any
$S \in \mathrm{GL}(\p(m, 1))$ can be extended to an automorphism $R(S)$ of $\f(m, p)$, any choice of an inner
product on $\p(m, 1)$ defines a nilsoliton inner product (if it exists) uniquely.

(a) Any derivation $\psi$ of $\f(m,p)$ which commutes with $\phi$ can be extended to a derivation $\hat \psi$ of the
metric solvable standard Einstein Lie algebra $\g =\mathbb{R}H \oplus \f(m,p)$ (orthogonal sum) by setting
$\hat \psi (H) = 0$. Then by \cite[Theorem 4.12]{H}, the endomorphism $\psi^*$ is also a derivation of $\f(m,p)$.
Taking $\psi = \rho(L)$ for some $L \in \End(\p(m, 1))$ we get two derivations, $(\rho(L))^*$ and
$\rho(L^*)$, which coincide on the generators, and hence are equal.

(b) This follows from assertion 1, if we take an endomorphism $L$ to be diagonal, $Le_i = \lambda_i e_i$, with the
$\la_i$'s linearly independent over $\mathbb{Q}$.

(c) This follows from the fact that the choice of the inner product on the generators determines the nilsoliton
inner product on $\f(m,p)$ uniquely (if the latter exists).

2. Let $\<\cdot,\cdot\>$ be an inner product on $\f(m, k)$, for which the spaces $\p(m, k)$ are orthogonal and
which satisfies (a), (b) and (c). Polarizing \eqref{eq:Riccinil} we obtain, for any $X \perp Y$,
\begin{equation}\label{eq:Ricpolarized}
\Ric_{\n} (X, Y) = \frac14 \sum\nolimits_{i,j} \<X, [E_i, E_j]\> \<Y, [E_i, E_j]\>
- \frac12 \sum\nolimits_{i} \< [X, E_i], [Y, E_i]\>,
\end{equation}
where $\{E_i\}$ is an orthonormal basis for $(\f(m, k), \<\cdot,\cdot\>)$.
The fact that the $\p(m, k)$'s are mutually orthogonal with respect to the quadratic form $\Ric$ follows immediately,
if we choose $\{E_i\}$ from $\cup_k \p(m, k)$. Next, in each of the $\p(m, k)$'s consider the \emph{content spaces}
$V_c$, the linear spans of the vectors $e_{i_1\ldots i_k}$ having the same content $c$. The content spaces $V_c$ are
mutually
orthogonal by (b) and their direct sum is the whole $\p(m, k)$ (as $\p(m, k)$ has a basis of the vectors of the form
$e_{i_1\ldots i_k}$). Choose an orthonormal basis $E_i$ in the union of the $L_c$'s. Then \eqref{eq:Ricpolarized}
and the fact that $[V_{c_1}, V_{c_2}] \subset V_{c_1+c_2}$ imply that $\Ric(V_{c_1}, V_{c_2}) =0$ when $c_1 \ne c_2$.

Finally, for a permutation $\tau$ of the generators $\{e_i\}$, let $L_\tau$ be the induced endomorphism of $\p(m,1)$.
Then, as it follows from (c), the automorphism $R(L_\tau)$ is an isometry, hence it preserves $\Ric$.
\end{proof}

In the following lemma, we explicitly compute the inner products on the spaces of Lie polynomials $\p(m, k)$ of those
free algebras which potentially may be Einstein nilradicals by Lemma~\ref{l:nonexist} (except for $\p(2, 7)$).
When $m=2$, we use the following notations: $\mathbf{e_k} = e_{121\ldots1} \in \p(2,k)$,
$\iota$ is an involution of $\p(2,k)$ defined by interchanging $1$ and $2$ in the subscripts
(that is, $\iota = R(S) \in \mathrm{Aut}(\f(2, p))$, where $S \in \mathrm{GL}(\p(2, 1))$ is defined by
$Se_1 = e_2,\; Se_2 = e_1$), $\Theta = \rho(L)$, where $L \in \End(\p(2,1))$ is defined by $Le_1 =e_2,\; Le_2 = 0$.

\begin{lemma}\label{l:product}
Any inner product $\<\cdot,\cdot\>$ on the space $\p(m, k)$ satisfying assertions (a), (b) and (c) of
Lemma~\ref{l:productgen} is uniquely determined by the choice of the corresponding nonzero
constants, as in the table below. The inner products of the basis vectors which are not listed can be either obtained
from those in the table by assertion (c) of Lemma~\ref{l:productgen}, or otherwise are zeros.

\begin{table}[h]
\setlength{\extrarowheight}{2pt}
\begin{center}
\begin{tabular}{|>{\cn}m{2.8cm}|>{\cn}m{1.8cm}|m{4cm}|m{5cm}|}
  \hline
   $\p(m, k)$ & $d_k(m)$ & \cn Basis & \parbox{4cm}{\cn Inner product} \\
  \hline
  $\p(m, 1)$ & $m$ & \cn $e_i$ & $\|e_{1}\| = 1$ \\
  \hline
  $\p(m, 2)$ & $\frac12(m^2-m)$ & \quad $e_{ij}, \hfill i < j$& $\|e_{12}\|^2 = \lambda^2$ \\
  \hline
  $\p(2, 3)$ & $2$ & \cn $e_{121}, e_{122}$ & $\|e_{121}\|^2 = \xi^2$ \\
  \hline
  $\p(m, 3),\; m > 2$ & $\frac13(m^3-m)$ & $e_{iji}, \hfill i \ne j$, $e_{ijk}, e_{jki}, \hfill i<j<k$ &
  $\|e_{121}\|^2 = \xi^2$, \newline
  $\|e_{123}\|^2 = \frac23 \xi^2$, \newline
  $\<e_{123}, e_{231}\> = -\frac13 \xi^2$\\
  \hline
  $\p(2, 4)$ &
  $3$ & \cn $e_{1211}, e_{1222}$, $e_{1212}$&
  $\|e_{1211}\|^2 = \sigma^2$, \newline
  $\|e_{1212}\|^2 = \frac12 \sigma^2$ \\
  \hline
  $\p(3, 4)$ &
  $18$ & $e_{ijii}, \hfill i \ne j$, \newline $e_{ijij}, \hfill i<j$, \newline
  $q_{ijk}= e_{kjii} + e_{ijki} + e_{ijik}$, \newline
  $q_{ikj}, \; r_{ijk}=[e_{ij},e_{ik}]$,
  \newline $i \ne j, k, \; j<k$ &
  $\|e_{1211}\|^2 = \sigma^2$, \newline
  $\|e_{1212}\|^2 = \frac12 \sigma^2$, \newline
  $\|q_{123}\|^2 = 3\sigma^2$, \newline
  $\|r_{123}\|^2 = \eta^2$, \newline
  $\<q_{123}, q_{132}\> = -\sigma^2$\\
  \hline
  $\p(2, 5)$ &
  $6$ & \cn $\mathbf{e_5}, \iota \mathbf{e_5}, \Theta \mathbf{e_5}, \iota \Theta \mathbf{e_5}$ \newline
  $u=[\mathbf{e_3},\mathbf{e_2}], \iota u$&
  $\|\mathbf{e_5}\|^2 = \alpha^2, \|\Theta \mathbf{e_5}\|^2 = 3\alpha^2$, \newline
  $\|u\|^2 = \gamma^2$ \\
  \hline
  $\p(2, 6)$ &
  $9$ &  \cn  $\mathbf{e_6}, \iota \mathbf{e_6}, \Theta \mathbf{e_6}, \iota \Theta \mathbf{e_6}, \Theta ^2\mathbf{e_6}$
  \newline
  $z=[\mathbf{e_4},\mathbf{e_2}], \; \Theta z, \; \iota z$,
  \newline $I = [\mathbf{e_3}, \Theta \mathbf{e_3}]$ &
  $\|\mathbf{e_6}\|^2 = \kappa^2, \|\Theta \mathbf{e_6}\|^2 = 4\kappa^2$,
  \newline $\|\Theta ^2\mathbf{e_6}\|^2 = 24\kappa^2$, \newline
  $\|z\|^2 = \delta^2, \|\Theta z\|^2 = 2 \delta^2$, \newline $\|I\|^2 = \theta^2$\\
  \hline
\end{tabular}
\caption{Inner product on the spaces $\p(m, k)$.}\label{table2}
\end{center}
\end{table}

\end{lemma}

\begin{proof}
The fact that the inner products on $\p(m, 2)$ and $\p(2, 3)$ have the form given in Table 2 immediately follows from
assertions (b) and (c) of Lemma~\ref{l:productgen}.

For $\p(m, 3)$ with $m \ge 3$, let $\|e_{121}\|^2 = \xi^2$ for some nonzero $\xi$. Define an endomorphism
$L$ of $\p(m, 1)$ by $Le_1 = e_3,\; Le_i = 0$ for $i \ne 1$. By assertions (a) of Lemma~\ref{l:productgen},
$\xi^2 = \<e_{121}, \rho(L^*)e_{123}\> = \<\rho(L)e_{121}, e_{123}\> = \<e_{321} + e_{123}, e_{123}\>$. Then from
$\|e_{123}\|^2 = \|e_{231}\|^2 = \|e_{312}\|^2$ (assertion (c) of Lemma~\ref{l:productgen}) and the Jacobi identity
$e_{123} + e_{231} + e_{312} = 0$ it follows that
$\|e_{123}\|^2 = \tfrac23 \xi^2, \; \<e_{123}, e_{231}\> = -\tfrac13 \xi^2$. All the other inner products of
the basis vectors can be then found from assertions (b) and (c) of Lemma~\ref{l:productgen}.

For $\p(2, 4)$, let $\|\mathbf{e_4}\|^2 = \sigma^2$. Then
$\sigma^2=\<\mathbf{e_4}, \Theta^*e_{1212}\> = \<\Theta \mathbf{e_4}, e_{1212}\>
= \<e_{1212} + e_{1221}, e_{1212}\> = 2\|e_{1212}\|^2$,
as $e_{1212}=e_{1221}$. The basis is orthogonal by assertion (b) of Lemma~\ref{l:productgen}.

For $\p(3, 4)$, take $\|e_{1211}\|^2 = \sigma^2$. The arguments similar to those for $\p(2, 4)$ show
that $2\|e_{1212}\|^2 = \sigma^2$. Define an endomorphism
$L$ of $\p(3, 1)$ by $Le_1 = e_3,\; Le_i = 0$ for $i \ne 1$. Then $q_{123}=\rho(L) e_{1211}$, so
$\|q_{123}\|^2=\<\rho(L^*) \rho(L) e_{1211}, e_{1211} \> = 3 \sigma^2$. For $r_{123}=[e_{12},e_{13}]$, set
$\|r_{123}\|^2 = \eta^2$. Then the norms of all the basis vectors are known by assertion (c)
of Lemma~\ref{l:productgen}. By assertions (b) and (c), the only inner products it remains to find are
$\<q_{123}, q_{132}\>$ and  $\<q_{123}, r_{123}\>$. The latter one is zero, as $q_{123}=\rho(L) e_{1211}$
and $\rho(L)^* r_{123} = 0$. For the former one, we have
$\<q_{123}, q_{132}\>= \<e_{1211}, \rho(L^*)(e_{2311} + e_{1321} + e_{1312})\> =
\<e_{1211}, e_{2111}\> = -\sigma^2$.

For $\p(2, 5)$, take $\|\mathbf{e_5}\|^2 = \alpha^2$. Then
$\|\Theta \mathbf{e_5}\|^2 = \<\Theta^*\Theta \mathbf{e_5}, \mathbf{e_5}\> = 3\|\mathbf{e_5}\|^2 = 3\alpha^2$. Denote
$u = [\mathbf{e_3},\mathbf{e_2}]$ and let
$\|u\|^2 = \gamma^2$. Then the norms of all the other basis vectors are determined by assertion (c)
of Lemma~\ref{l:productgen}, and by assertion (b), any two of them are orthogonal, except possibly $A\mathbf{e_5}$ and
$u$ (and $\iota \Theta \mathbf{e_5}$ and $\iota u$). However, as $\Theta^*u = 0, \quad \<u, \Theta \mathbf{e_5}\> = 0$.

For $\p(2, 6)$, take $\|\mathbf{e_6}\|^2 = \kappa^2$. Then
$\|\Theta \mathbf{e_6}\|^2 = \<\Theta^*\Theta \mathbf{e_6}, \mathbf{e_6}\> = 4\|\mathbf{e_6}\|^2 = 4\kappa^2$, and by
a similar computation $\|\Theta ^2\mathbf{e_6}\|^2 = 24\kappa^2$. Denote $z = [\mathbf{e_4},\mathbf{e_2}]$ and set
$\|z\|^2 = \delta^2$. As $\Theta^*\Theta z = 2z$, we obtain $\|\Theta z\|^2 = 2\kappa^2$. Denote
$I = [\mathbf{e_3}, \Theta \mathbf{e_3}]= [e_{121}, e_{122}]$ and set $\|I\|^2 = \theta^2$.
The norms of all the basis vectors are now known (assertion (c) of Lemma~\ref{l:productgen}). Also, the
subspaces $\Span(\mathbf{e_6}), \, \Span(\Theta \mathbf{e_6}, z), \, \Span(\Theta ^2\mathbf{e_6}, \Theta z, I)$,
$\Span(\iota \Theta \mathbf{e_6}, \iota z)$, and $\Span( \iota \mathbf{e_6})$ are orthogonal by assertion (b)
of Lemma~\ref{l:productgen}. As $\Theta^* z = 0, \; z \perp \Theta \mathbf{e_6}$ (and
$\iota z \perp \iota \Theta \mathbf{e_6}$).
Similarly, $\Theta^* I = 0$, so $I \perp \Theta ^2\mathbf{e_6}, \Theta z$. From $\Theta^*\Theta z = 2z$ it follows that
$(\Theta^2)^*\Theta z = 0$ which shows that $\Theta ^2\mathbf{e_6} \perp \Theta z$.

\end{proof}

In order for $\f(m, p)$ to be an Einstein nilradical, the inner product on it constructed according to
Lemma~\ref{l:productgen} and Lemma~\ref{l:product} has to satisfy \eqref{eq:ricn}. By equation \eqref{eq:Riccinil} and
Lemma~\ref{l:der}, this is equivalent to the fact that for all $X \in \p(m, k)$,
\begin{equation}\label{eq:ricXX}
\Ric_{\n} (X) = \frac14 \sum_{i,j} \<X, [E_i, E_j]\>^2 - \frac12 \sum_{i} \| [X, E_i]\|^2 =
(k \, \Tr \, \Phi - \Tr \, \Phi^2) \, C \, \|X\|^2,
\end{equation}
where $C = -c (\Tr \, \Phi^2)^{-1} > 0$ (recall that $\Phi$ is the canonical derivation defined by \eqref{eq:canPhi},
$\Tr \, \Phi = \sum_k k d_k(m)$, $\Tr \, \Phi^2 =  \sum_k k^2 d_k(m)$, where the $d_k(m)$'s are given in
Table~\ref{table}).

\bigskip

To finish the proof of the Theorem, it remains to compute the Ricci tensor of those free Lie algebras which can
potentially be Einstein nilradicals by Lemma~\ref{l:nonexist} and to check, whether the
constants in Table~\ref{table2} can be chosen in such a way that \eqref{eq:ricXX} is satisfied.

It is well-known that any abelian Lie algebra is an Einstein nilradical
(\cite[Proposition 4.2]{A2}, \cite[Proposition 6.12]{H}). The fact that any two-step free Lie algebra is an Einstein
nilradical easily follows from \eqref{eq:ricXX} (or from \cite[Proposition 2.9(iii)]{GK}).

In all the other cases, we chose the inner product according to Lemma~\ref{l:product} and normalized by the condition
$\|e_{12}\| = 1$.

\bigskip

\noindent $\mathbf{\f(m, 3)}$

Equation \eqref{eq:ricXX} gives
\begin{align*}
\Ric_{\n} (e_1) &= - \tfrac12(m-1) - \tfrac{1}{6} \xi^2 (m^2-1) = (- 2m^3 -m^2 +3m) C, \\
\Ric_{\n} (e_{12}) &= \tfrac12 - \tfrac{1}{3} \xi^2 (m+1) = (-m^3 + 2m) C, \\
\Ric_{\n} (\xi^{-1} e_{121}) &= \Ric_{\n} (\sqrt{\tfrac32} \xi^{-1} e_{123}) = \tfrac 12 \xi^2 = (m^2 + m) C,
\end{align*}
from which
$\xi^2 = 3(m+1)(-m^2 + 4m + 8)^{-1}$, so $\f(m, 3)$ can be an Einstein nilradical only when $-m^2 + 4m + 8 > 0$,
(this inequality already appeared in the proof of Lemma~\ref{l:nonexist}),
that is, when $m = 2, 3, 4, 5$. In all these
cases, $\f(m, 3)$ is indeed an Einstein nilradical by assertion 2 of Lemma~\ref{l:productgen}.

%\bigskip

\newpage

\noindent $\mathbf{\f(2, 4)}$

From equation \eqref{eq:ricXX} we get (recall that $\mathbf{e_k} = e_{121\ldots1}=[\ldots[e_1,e_2],e_1],\ldots],e_1]$,
$k-1$ ones in total):
\begin{align*}
\Ric_{\n} (e_1) &= - \tfrac12 - \tfrac12 \xi^2 - \tfrac34 \xi^{-2} \sigma^2(m-1) = -50 C, \\
\Ric_{\n} (\mathbf{e_2}) &= \tfrac12 - \xi^2 = -28 C, \\
\Ric_{\n} (\xi^{-1} \mathbf{e_3}) &= \tfrac12 \xi^2 - \tfrac34 \xi^{-2} \sigma^2 = -6 C,\\
\Ric_{\n} (\sigma^{-1} \mathbf{e_4}) = \Ric_{\n} (\sqrt2 \sigma^{-1} e_{1212}) &= \tfrac12 \xi^{-2} \sigma^2 = 16 C,
\end{align*}
which gives $\xi^2 = \frac94,\; \sigma^2 = \frac92$, so $\f(2, 4)$ is an Einstein nilradical.

\bigskip

\noindent $\mathbf{\f(3, 4)}$

Equation \eqref{eq:ricXX} gives
\begin{align*}
\Ric_{\n} (e_{12}) &= \tfrac12 - \eta^2 - \tfrac43 \xi^2 = -165 C, \\
\Ric_{\n} (\xi^{-1} e_{121}) &=
\tfrac12 \xi^2 - \tfrac{15}{16} \xi^{-2} \sigma^2 - \tfrac{9}{32} \xi^{-2} \eta^2= -60 C,\\
\Ric_{\n} (\sigma^{-1} e_{1211}) &= \tfrac12 \xi^{-2} \sigma^2 = 45 C,\\
\Ric_{\n} (\eta^{-1} [e_{12},e_{13}]) &= \tfrac12 \eta^2 + \tfrac34 \xi^{-2} \eta^2 = 45 C,
\end{align*}
from which $16 \xi^4 + 171 \xi^2 + 99 = 0$, a contradiction.

\bigskip

\noindent $\mathbf{\f(2, 5)}$

From equation \eqref{eq:ricXX} we get
\begin{align*}
\Ric_{\n}(e_1) &= -\tfrac12 - \tfrac12 \xi^2 - \tfrac34 \xi^{-2} \sigma^2
- \sigma^{-2} (\alpha^2 + \tfrac{1}{3} \gamma^2) = -170 C,\\
\Ric_{\n}(\mathbf{e_2}) &= \tfrac12 - \xi^2 - \xi^{-2} \gamma^2 = -118 C,\\
\Ric_{\n}(\xi^{-1}\mathbf{e_3}) &= \tfrac12 \xi^2 - \tfrac34 \xi^{-2} \sigma^2 - \tfrac12 \xi^{-2} \gamma^2 = -66 C,\\
\Ric_{\n}(\sigma^{-1}\mathbf{e_4}) = \Ric_{\n} (\sqrt2 \sigma^{-1} e_{1212})
&= \tfrac12 \xi^{-2} \sigma^2 - \tfrac{2}{9} \sigma^{-2} \gamma^2 - \tfrac{2}{3} \sigma^{-2} \alpha^2 = -14 C, \\
\Ric_{\n}(\gamma^{-1}u) &= \tfrac{1}{3} \sigma^{-2} \gamma^2 + \tfrac12 \xi^{-2} \gamma^2  = 38 C, \\
\Ric_{\n}(\alpha^{-1}\mathbf{e_5}) = \Ric_{\n}((\sqrt3 \alpha)^{-1} \Theta \mathbf{e_5})
&= \tfrac12 \sigma^{-2} \alpha^2 = 38 C. \\
\end{align*}
Solving this we obtain $\xi^2 = 54 C, \; \sigma^2 = 3/4+375 C, \; \alpha^2 = 76 C \sigma^2, \;
\gamma^2 = 27 C (1 + 128 C)$, where $C=(131+5\sqrt{745})/2928 \approx 0.09135$ is the positive root of
$5856 C^2-524 C-1=0$.

\bigskip

\noindent $\mathbf{\f(2, 6)}$

By \eqref{eq:ricXX} we get
\begin{equation}\label{eq:ricf(2,6)}
\begin{split}
\Ric_{\n}(\mathbf{e_2}) &= \tfrac12 - \xi^2 - \xi^{-2} \gamma^2 - \tfrac32 \sigma^{-2} \delta^2= -334 C\\
\Ric_{\n}(\xi^{-1}\mathbf{e_3}) &= \tfrac12 \xi^2 - \tfrac34 \xi^{-2} \sigma^2 - \tfrac12 \xi^{-2} \gamma^2
- \tfrac12 \xi^{-4} \theta^2 =-228 C\\
\Ric_{\n}(\sigma^{-1}\mathbf{e_4}) &= \tfrac12 \xi^{-2} \sigma^2 - \tfrac{2}{9} \sigma^{-2} \gamma^2
- \tfrac{2}{3} \sigma^{-2} \alpha^2 - \tfrac12 \sigma^{-2} \delta^2 = -122 C \\
\Ric_{\n}(\gamma^{-1}u) &= \tfrac{1}{3} \sigma^{-2} \gamma^2 + \tfrac12 \xi^{-2} \gamma^2
- \tfrac12 \gamma^{-2} \theta^2 - \tfrac34 \gamma^{-2} \delta^2 = -16 C \\
\Ric_{\n}(\alpha^{-1}\mathbf{e_5}) &= \tfrac12 \sigma^{-2} \alpha^2 - \tfrac{5}{8}\alpha^{-2} \kappa^2
 - \tfrac{25}{32}\alpha^{-2} \delta^2= -16 C \\
\Ric_{\n}(\theta^{-1}I) &= \gamma^{-2} \theta^2 + \tfrac12 \xi^{-4} \theta^2= 90 C \\
\Ric_{\n}(\delta^{-1}z) &= \tfrac12 \gamma^{-2} \delta^2 + \tfrac12 \sigma^{-2} \delta^2
+ \tfrac{25}{24}\alpha^{-2} \delta^2= 90 C \\
\Ric_{\n}(\kappa^{-1}\mathbf{e_6}) &= \tfrac12 \alpha^{-2} \kappa^2 = 90 C. \\
\end{split}
\end{equation}
This leads to a contradiction, as
\begin{multline*}
2 \, \Ric_{\n}(\xi^{-1}\mathbf{e_3}) + 3 \, \Ric_{\n}(\sigma^{-1}\mathbf{e_4}) +
2 \, \Ric_{\n}(\gamma^{-1}u) + 4 \, \Ric_{\n}(\alpha^{-1}\mathbf{e_5}) \\
+ 2 \, \Ric_{\n}(\theta^{-1}I) + 3 \, \Ric_{\n}(\delta^{-1}z)+ 5 \, \Ric_{\n}(\kappa^{-1}\mathbf{e_6}) =
\xi^2 + \gamma^{-2} \theta^2 = -18 C <0.
\end{multline*}

\bigskip

\noindent $\mathbf{\f(2, 7)}$

We start by computing the inner product on $\p(2, 7)$.
The space $\p(2, 7)$ has dimension $18$. Choose a basis for $\p(2, 7)$ consisting of the nine vectors
$$
\mathbf{e_7}, \Theta \mathbf{e_7}, [\mathbf{e_3}, \mathbf{e_4}], [\mathbf{e_2}, \mathbf{e_5}],
\Theta ^2\mathbf{e_7}, \Theta [\mathbf{e_3}, \mathbf{e_4}], \Theta [\mathbf{e_2}, \mathbf{e_5}],
[u,\mathbf{e_2}], [I, e_1]
$$
and the nine vectors obtained from them by the action of $\iota$. From assertion (b) of Lemma~\ref{l:productgen} it
follows that the six subspaces $\mathbb{R}\mathbf{e_7}$,
$\Span(\Theta \mathbf{e_7}, [\mathbf{e_3}, \mathbf{e_4}], [\mathbf{e_2}, \mathbf{e_5}])$,
$\Span(\Theta ^2\mathbf{e_7}$, $\Theta [\mathbf{e_3}, \mathbf{e_4}],
\Theta [\mathbf{e_2}, \mathbf{e_5}], [u,\mathbf{e_2}], [I, e_1])$ and the other three obtained from
these by the action of $\iota$ are mutually orthogonal.

Set $\|\mathbf{e_7}\|^2 = \nu^2$. From assertion (a) of Lemma~\ref{l:productgen}, $\|\Theta \mathbf{e_7}\|^2 = 5 \nu^2$
and $\|\Theta^2 \mathbf{e_7}\|^2 = 40 \nu^2$. Let
$V = \Span([\mathbf{e_3}, \mathbf{e_4}], [\mathbf{e_2}, \mathbf{e_5}])$. As
$\Theta^* V = 0, \quad \Theta \mathbf{e_7} \perp V$. Furthermore, for any $X \in V, \quad \Theta^* \Theta X = 3 X$, so
$\|\Theta X\|^2 = 3 \|X\|^2$ and $\Theta X \perp \Theta^2 \mathbf{e_7}$. Let $W = \Span([u,\mathbf{e_2}], [I, e_1])$. As
$\Theta^* W = 0$, we get $\Theta^2 \mathbf{e_7}, \Theta V \perp W$.

Denote
$\Ric_{\n}(V)$ (respectively, $\Ric_{\n}(W)$) the trace of the restriction of the quadratic form $\Ric$ to $V$
(respectively, to $W$). A direct computation using the inner products given in Table~\ref{table2} and
equation \eqref{eq:ricXX} gives
\begin{align*}
\Ric_{\n}(V) &= \tfrac12 \xi^{-2} \sigma^{-2} \|[\mathbf{e_3}, \mathbf{e_4}]\|^2
+ \tfrac12 \alpha^{-2} \|[\mathbf{e_2}, \mathbf{e_5}]\|^2 +
\tfrac54 \kappa^{-2} \|[\mathbf{e_6}, e_2]_V\|^2
+ \tfrac12 \delta^{-2} \|[z, e_1]\|^2 = 392 C, \\
\Ric_{\n}(W) &= 3 \xi^{-2} \sigma^{-2} \|[\mathbf{e_3}, e_{1212}]_W\|^2 + \tfrac12 \gamma^{-2} \|[\mathbf{e_2}, u]\|^2 +
\tfrac12 \theta^{-2} \|[I, e_1]\|^2
+ \tfrac34 \delta^{-2} \|[z, e_2]_W\|^2  = 392 C,\\
\Ric_{\n}(\nu^{-1} \mathbf{e_7}) &= \tfrac12 \kappa^{-2} \nu^2 = 196 C,\\
\Ric_{\n}(\mathbf{e_2}) &= \Ric_{\f(2, 6)}(\mathbf{e_2}) - 2 \alpha^{-2} \|[\mathbf{e_2}, \mathbf{e_5}]\|^2
- \gamma^{-2} \|[\mathbf{e_2}, u]\|^2 = -964 C\\
\Ric_{\n}(\xi^{-1}\mathbf{e_3}) &= \Ric_{\f(2, 6)}(\xi^{-1}\mathbf{e_3})
- \xi^{-2} \sigma^{-2} (\|[\mathbf{e_3}, \mathbf{e_4}]\|^2 + 3 \|[\mathbf{e_3}, e_{1212}]_W\|^2) =-732 C\\
\Ric_{\n}(\sigma^{-1}\mathbf{e_4}) &= \Ric_{\f(2, 6)}(\sigma^{-1}\mathbf{e_4})
- \xi^{-2} \sigma^{-2} (\tfrac23 \|[\mathbf{e_3}, \mathbf{e_4}]\|^2
+ \|[\mathbf{e_3}, e_{1212}]_W\|^2) = -500 C \\
\Ric_{\n}(\gamma^{-1}u) &= \Ric_{\f(2, 6)}(\gamma^{-1}u) - \tfrac12 \gamma^{-2} \|[\mathbf{e_2}, u]\|^2 = -268 C \\
\Ric_{\n}(\alpha^{-1}\mathbf{e_5}) &= \Ric_{\f(2, 6)}(\alpha^{-1}\mathbf{e_5})
- \tfrac12 \alpha^{-2} \|[\mathbf{e_2}, \mathbf{e_5}]\|^2 = -268 C \\
\Ric_{\n}(\theta^{-1}I) &= \Ric_{\f(2, 6)}(\theta^{-1}I) - \theta^{-2} \|[I, e_1]\|^2 = -36 C \\
\Ric_{\n}(\delta^{-1}z) &= \Ric_{\f(2, 6)}(\delta^{-1}z) - \delta^{-2}( \tfrac23 \|[z, e_1]\|^2
+ \tfrac12 \|[z, e_2]_W\|^2) = -36 C \\
\Ric_{\n}(\kappa^{-1}\mathbf{e_6}) &= \Ric_{\f(2, 6)}(\kappa^{-1}\mathbf{e_6})
- \kappa^{-2} (\tfrac35 \nu^2 + \tfrac12 \|[\mathbf{e_6}, e_2]_V\|^2) = -36 C,
\end{align*}
where $\Ric_{\f(2, 6)}$ are the corresponding expressions from \eqref{eq:ricf(2,6)} and the subscripts $V$
and $W$ mean the orthogonal projections to the corresponding subspace. This implies
\begin{multline*}
\Ric_{\n}(\mathbf{e_2}) + 2  \Ric_{\n}(\xi^{-1}\mathbf{e_3}) + 6  \Ric_{\n}(\sigma^{-1}\mathbf{e_4}) +
8  \Ric_{\n}(\alpha^{-1}\mathbf{e_5}) + 4  \Ric_{\n}(\gamma^{-1}u) \\
+ 10  \Ric_{\n}(\kappa^{-1}\mathbf{e_6}) +
9  \Ric_{\n}(\delta^{-1}z) + 3  \Ric_{\n}(\theta^{-1}I) +
12  \Ric_{\n}(V) + 6  \Ric_{\n}(W) + 12  \Ric_{\n}(\nu^{-1} \mathbf{e_7}) \\
= \tfrac12 + \tfrac32 \sigma^2 \xi^{-2} + \gamma^{-2} \theta^2 + \tfrac12 \xi^{-4} \theta^2
+ \tfrac32 \gamma^{-2} \delta^2 +
\tfrac{25}{8} \alpha^{-2} \delta^2 + \tfrac52 \kappa^{-2} \|[\mathbf{e_6}, e_2]_V\|^2 = -28 C,
\end{multline*}
which is a contradiction.

\end{document}